\documentstyle[a4,11pt,amssymb,leqno]{amsart}

\title[Metrically thin singularities of 
integrable CR functions]
{Metrically thin singularities of 
integrable CR functions}

\def\R{{\mathbb R}} 
\def\C{{\mathbb C}}

\date{\today}

\author{Jo\"el Merker and Egmont Porten}

\address{Laboratoire d'Analyse, Topologie et Probabilit\'es,
Centre de Math\'ematiques et d'Informatique, UMR 6632,
39 rue Joliot Curie,
F-13453 Marseille Cedex 13, France}

\email{merker@@cmi.univ-mrs.fr, merker@@dmi.ens.fr}

\address{ Max-Planck-Gesellschaft, Humboldt-Universit\"at zu Berlin 
J\"agerstrasse, 10-11, D-10117 Berlin, Germany} 

\email{egmont@@mathematik.hu-berlin.de}



\newtheorem{thm}{Theorem}[section]

\newtheorem{lem}[thm]{Lemma}
\newtheorem{rem}[thm]{Remark}

\newcommand{\crdim}{\mbox{CR-dim}}
\newcommand{\eps}{\epsilon}

\def\R{{\bf R}}
\def\C{{\bf C}}

\begin{document}

\begin{abstract} 
In this article, we consider metrically thin singularities $A$ 
of the tangential Cauchy-Riemann operator on smoothly embedded
Cauchy-Riemann manifolds $M$. The main result states removability
within the space of locally integrable functions on $M$ under
the hypothesis that the $(\dim M-2)$-dimensional Hausdorff volume
of $A$ is zero and that the CR-orbits of $M$ and $M\backslash A$ 
are comparable. 
\end{abstract}

\maketitle

\section{Introduction}

A smooth real submanifold $M$ of a complex manifold $X$
is called (embedded) CR-manifold
if the dimension of the maximal complex subspace $T^c_p M$
of the tangent space $T_p M$ does not depend on $p\in M$.
In this case the complex tangent spaces are the fibers of
a smooth vector bundle, the complex tangent bundle $T^c M$,
whose (complex) rank $\crdim M=\dim_c T^c_p M$ is called the
CR-dimension of $M$. Then the CR-vectorfields, i.e.~the
sections of $T^c M$, form a system of first order differential
operators, also denoted as the system of the tangential
Cauchy-Riemann equations. Functions annihilated
by all CR-vectorfields will be called CR-functions.

This paper is devoted to the study of the singularities
of integrable CR-functions. A closed set $A\subset M$
is called $L^{\rm p}$-removable if
\[
  L_{loc,CR}^{\rm p}(M\backslash A) \cap
  L_{loc}^{\rm p}(M)= L_{loc,CR}^{\rm p}(M),
\]
{\it i.e.} if any $f\in L_{loc}^{\rm p}(M)$
which is CR outside of $A$ is automatically
CR on all of $M$. This notion appeared probably
for the first time in the classical Riemann-Removability
Theorem (stating the $L^\infty_{loc}$ removability of isolated
points for the Cauchy-Riemann operator on the complex plane).
Later on, the question was treated by Bochner, Carleson
e.a.~also in the other $L^{\rm p}$-spaces for the ordinary
Cauchy-Riemann and the Laplace operator, and by
Harvey and Polking for general linear partial differential
operators. Their main theorem (Theorem 4.1.~of [HAPO]),
which is best possible in the general setting,
implies for tangential Cauchy-Riemann operators
that $E$ is $(L^{\rm p},\overline{\partial}_{\rm b})$-removable
if the Hausdorff measure $H^{\dim M-p'}(E)<\infty$, $p'=p/(p-1)$.
We stress that on the general level considered by Harvey and
Polking no results in $L^1$ are possible, as is already shown
by the meromorphic function of one variable $1/z$, which is
locally integrable near the singularity at the origin.

The main theorem of this article shows that essentially stronger
removability phenomena hold true for tangential Cauchy-Riemann
operators on embedded CR-manifolds which are not Riemann surfaces.
Our approach relies on
the close interplay between CR-functions and complex analysis
by means of analytic extension. The rough idea tells that
singularities of CR-functions should behave in a way comparable
to holomorphic functions. For boundaries of complex domains,
this point of view was illustrated by contributions
of many authors (we refer to the excellent survey [CS]).
In the present paper, we focus in particular on CR-manifolds
of arbitrary codimension.

As background for analytic extension we need to recall the
notion of CR-orbits. For a point $p$ of a CR-manifold
$M$ its CR-orbit ${\cal O}(p,M)$ is defined as the union of
all points which are connected with $p$ by a piecewise
differentiable path $\gamma\subset M$
whose (one-sided) derivatives are
always non-vanishing and contained in $T^c M$.
By general results of
Sussmann it is known that CR-orbits are injectively
immersed submanifolds of $M$. For information on their
tremendous bearing on analytic extension and further
references, we refer the reader to contributions of
Treves, Tr\'epreau, Tumanov, J\"oricke and Merker.

We shall constantly employ the following terminology:
A property $P$ is said to be true on almost every CR-orbit
if there is a union $\cal N$ of CR-orbits such that
$P$ holds true on any orbit $\cal O$ out of $\cal N$.
Furthermore a CR-manifold $M$ will be called globally minimal
if it contains only one orbit.

The following theorem contains the essence of the
article.

\begin{thm}\label{main}
Let $M$ be an embedded CR-manifold of class
${\cal C}^{3}$,
%, $0< \alpha <1$,
dimension $d\geq 3$, and CR-dimension
${\rm dim}_{CR} M =m\geq 1$.
Then every closed subset $E$ of $M$ with
$H^{d-2}(E)=0$ such that for almost all
CR orbits, ${\cal O}_{CR}\backslash E$ is globally minimal
is $L^p_{loc}$-removable for $1\leq p\leq\infty$.
\end{thm}

The above theorem is a refinement of results the authors
proved in [MP], where either the
hypothesis on Hausdorff measure was too restrictive
(we supposed $H^{d-3}(E)$ locally finite), or $M$ was
assumed to be real analytic. This research was as a whole
inspired by corresponding results of Lupacciolu, Stout and Chirka
on analytic extension of continuous CR-functions defined on
parts of boundaries of complex domains.

The metrical hypothesis is optimal for trivial reasons:
Taking for $E$ the transverse intersection of $M$ with with
a complex hypersurface $X$ of the ambient space we plainly
get $H^{d-2}(E)>0$. If $X$ has a holomorphic defining
function $f$, then $1/f|_M$ is around $E$ locally integrable  
but not CR.

The organization of the proof is as follows: In a first
reduction (section 2) we shall see that the problem restricts
to solving the
corresponding question on almost all CR-orbits. More
concretely we shall be left with the hypothesis that
both $M$ and $M\backslash (M-E)$ are globally minimal.

Next we employ known techniques using Bishop-discs to
extend our CR-function $f$ from $M\backslash E$ to a large
portion of a wedge $W$ attached to $M$ (section 3).
Then the final step should be to remove the singularity
in the wedge with good $L^1$-estimates and to recover
thereafter $f$ globally as an $L^1$-limit (similary as in the
theory of Hardy spaces) thereby proving that $f$ is
CR everywhere. At the removal of the singularity in $W$
we encounter a special problem which is overcome by
the following theorem which seems to be of independent
interest.

\begin{thm}\label{hol}
Let $D\subset\C^N$ be a domain
equiped with a foliation $\cal F$ by holomorphic curves
of class $C^2$. Further let $E\subset D$ be a closed union
of leaves with $H_{2N}(E)=0$. Then a function
$F\in L^1(D) \cap {\cal H}(D\backslash E)$
extends holomorphically to $D$ as soon as $F$ extends
to an open subset $D'$ of $D$ which intersects each
leaf contained in $E$.
\end{thm}

The proof is an extension of techniques used by Henkin and
Tumanov to prove a related statement on continuous CR-functions
on manifolds foliated by holomorphic curves.

% If $\cal F$ in Theorem \ref{hol} is even a {\sl holomophic} foliation,
% there is an elementary proof by the continuity principle,
% which moreover does not use the integrability of $F$.
% For dimensional reasons
% The seemingly
% special case $N=2$ generalizes readily to differential
% foliations by analytic {\sl hypersurfaces}. For the general situation
% of the above theorem we do not know if a proof working only
% with holomorphic hulls exists.

We would like to thank B.~J\"oricke for valuable discussions.

\section{Reduction to CR-orbits}

In this section we shall reduce the proof of theorem \ref{main}
to the case that both $M$ and $M\backslash E$ are globally minimal.
This is a consequence of the following two statements.

\begin{lem}\label{l1}
Under the assumptions of Theorem \ref{main}, for almost every
orbit $\cal O$ the intersection $E\cap{\cal O}$ satisfies
$H_E^{d_{\cal O}-2}(E\cap{\cal O})=0$, where 
$d_{\cal O}= {\rm dim}_{\R} {\cal O}.$
\end{lem}

Here the index $E$ in $H_E^{d_{\cal O}-2}(E\cap{\cal O})$ indicates
that the Hausdorff content is computed with respect
to the manifold topology of $\cal O$
and the pullback of the euclidean metric to $\cal O$.

{\it Proof.}
First we recall the following local property of CR-orbits:
Near every point $p\in M$ there exist a (topologically trivial)
$C^2$-foliation of $M$
into  leaves of the dimension $k=\dim{\cal O}(p,M)$ such that
for any $q$ near $p$ its orbit ${\cal O}(q,M)$ contains
the leaf passing through $q$.
This implies that the dimension
of the CR-orbits passing near $p$ cannot
be less than $\dim{\cal O}(p,M)$ and that in case of equality
the leaf is an open subset of the corresponding orbit.
A well-known property of Hausdorff-measures assures that
for almost every leaf $L$ the intersection $E\cap L$
satisfies $H^{\dim L-2}(E\cap L)=0$.

To exploit this local information we need the
following simple covering lemma.

\begin{lem}\label{l11}
There is a countable covering of $M$ by
foliated open sets $U_j\subset\subset M$
as above such that every point $p\in M$
appears in at least one $U_j$ whose leaves are of dimension
$\dim{\cal O}(p,M)$.
\end{lem}

{\it Proof.}
We argue by induction on the dimension of orbits.
By paracompactness of $M$ we easily find a countable covering
$\{U_j^{(0)}\}$ as required for the open submanifold
\[
  M^{(0)}=\{ p\in M:\dim{\cal O}(p,M)=\dim M\}.
\]
In this case the $U_j^{(0)}$  are just equiped with the
trivial codimension-zero-foliation.

In the next we look at
\[
  M^{(1)}=\{ p\in M:\dim{\cal O}(p,M)\geq\dim M-1\},
\]
which is by the above local property of CR-orbits also
an open subset of $M$. We choose for any point
$p\in M_1\backslash M_0$ an open set $U_p^{(1)}$ equiped
with a convenient codimension-one-foliation and obtain
thereby a covering $\{U_j^{(0)},U_p^{(1)}\}$. Again by
paracompactness we extract a countable covering
$\{U_j^{(1)}\}\supset\{U_j^{(0)}\}$.

By induction we treat the open submanifolds
\[
  M^{(k)}=\{ p\in M:\dim{\cal O}(p,M)\geq\dim M-k\}
\]
in analogous fashion. Hence we get an increasing
chain of systems $\{U_j^{(k)}\}$ and the desired covering
after at most $\dim M-2\crdim M$ steps. \hfill $\square$

From Lemma \ref{l11} and the preceding remark we deduce
that the set of all CR-orbits $\cal O$ of a given
dimension $k$ for which
the intersection ${\cal O}\cap E$ with some leaf of
with some fixed $U_j$ (with leaf-dimension $k$) 
has too large Hausdorff content is of zero
measure. As countable unions of sets of measure zero
are still of measure zero, Lemma \ref{l1} is proved. \qed

Next we need the following result from [PO].

\begin{thm}\label{t1}
If $M$ is ${\cal C}^3$,
a function $f\in L_{loc}^1(M)$ is CR if and only if $f|_{{\cal O}}$
belongs to $L_{loc,CR}^1({\cal O})$ for almost every CR-orbit
${\cal O}$.
\end{thm}

Lemma \ref{l1} obviously implies that the set-theoretical
assumptions of Theorem \ref{main} are inherited by almost
all orbits. Theorem \ref{t1} shows, first, that the restriction
of an $f\in L_{loc,CR}^1({\cal M\backslash E})$ 
is for almost every orbit a CR-function on
${\cal O}\backslash E$ and, secondly, that it is enough
to remove the singularity orbit-wise. Hence the desired reduction
to globally minimal manifolds is complete.

\section{Reduction to singularities of holomorphic functions}

Now we shall use the technique of analytic discs to reduce
the proof to a special problem concerning removable singularities
of integrable holomorphic functions. By the preceding section
and the fact that the CR-orbits of a manifold $M$ of class $C^{2,\alpha},
0<\alpha<1$, are of class $C^{2,\beta}$, for
any $\beta\in(0,1)$ (s.~[MP]),
it is enough to prove the following theorem.

\begin{thm}
Let $M$ be ${\cal
C}^{2,\alpha}$, $0< \alpha <1$, ${\rm dim}_{CR} M =m\geq 1$,
$\dim M\geq 3$. Then every closed subset $E$ of $M$
such that $M$ and $M\backslash E$ are globally
minimal and such that $H_{loc}^{d-2}(E)=0$, is $L^1_{loc}$-removable.
\end{thm}

{\it Proof.} As the proof is quite long, it will subdivided
into several steps:

{\it Step one: general setup.}
Fix a function $f\in L^1_{loc,CR}(M\backslash E)$.
Define
\[
{\cal A}= \{\Psi \subset E \ \hbox{closed} ; M\backslash \Psi
\ \hbox{globally minimal and} \ f\in L_{loc,CR}^1(M\backslash \Psi)
\cap L_{loc}^1(M)\}
\]
and define $E_{{\rm nr}}=\cap_{\Psi\in {\cal A}}
\Psi$.  Then $M\backslash E_{{\rm nr}}$ is globally minimal too.
If $E_{{\rm nr}}=\emptyset$, we are done.

To reach a contradiction, we assume $E_{{\rm nr}}\not=\emptyset$
and denote $E_{{\rm nr}}$ from now on again by $E$.
As $M\backslash E$ is globally minimal, we may, according to
Proposition 1.16 of
[MP1], (after a slight deformation of
$M$ fixing $E$) assume that
$f\in L_{loc}^1(M)\cap {\cal H}({\cal V}(M\backslash E))$.

First we look for a small analytic disc $A$ attached to $M$
whose boundary has non-trivial intersection with $E$ without
being contained in $E$. Further we shall achieve that $A$
is round in the following sense: There is a point $p\in bA$
and holomorphic coordinates $z,w$ around $p$ exhibiting
such that $\{z=0\}$ corresponds to $T^c_p M$ such that the
projection of $A$ to  $\{z=0\}$ is a translated coordinate disc.

Indeed, let $q\in E$ first be an arbitrary point
and $\gamma$ a piece-wise
differentiable CR-curve linking $q$ with a point
$q'\in M\backslash E$. After shortening $\gamma$ we may
suppose $\{q\}=E\cap\gamma$ and $\gamma$ to be a
smoothly embedded segment. Therefore $\gamma$ can
be described as the integral curve of some non-vanishing
CR vectorfield $X$ defined in a neighborhood of $\gamma$.

It suffices to show that $f$ is CR near
the endpoint $q$. Indeed, a standard argument using
the dynamical flow of $X$ proves the existence of
arbitrarily small neighborhoods $V$ of $q$ for which
$(M\backslash E)\cup V$ is globally minimal.
If $f$ is CR near $q$ and $V$ small enough, this
contradicts the minimality of $E$.

{\it Step two: construction of analytic discs.}
Let us briefly recall some basic notions concerning
analytic discs. An analytic disc $A$ is a holomorphic mapping
$A:\Delta\rightarrow\Omega$ extending continuously to
$\overline\Delta$. Sometimes, we will by $A,\partial A$ also denote the
images $A(\Delta),A(\partial\Delta)$.  The disc $A$ is said to be
attached to $M$ if $\partial A\subset M$.

Locally, we will usually work in holomorphic euclidean coordinates
$(z=x+iy,w)$ transfering a given point $p\in M$ to the origin such
that $M$ is given near $p$ as a graph of a function
\begin{equation}
  x = h(y,w),
\end{equation}
satisfying $h(0,0)=0$ and $\nabla h(0,0)=0$.  A disc $A$ will be
called round if in appropriated coordinates $(z,w)$ its $w$-components
appear as a round disc in some complex line in $\C^m$.

By an observation from [TU] we may approximate $\gamma$ by a finite
chain of round analytic discs
$A_1,\ldots,A_k,A_1(1)=q,A_j(-1)=A_{j+1}(1)$, $A_k(-1)$ very close to
$q'$, where all $A_j$ are of small diameter.

In what follows we shall prove for a round disc $A$ of small size
that $f$ is CR near all of $\partial A$ as soon as $\partial
A\not\subset E$.  If $A_j$ is the first of the above discs whose
boundary is not contained in $E$ this will imply that $f$ is CR near
$\partial A_j$. Careful examination of the constructions below yields that
we can successively apply the argument to $A_{j+1},\ldots,A_k$. In
particular $f$ is CR near $q$ leading to a contradiction.

{\it Step three: partial analytic extension.}  For a small round disc
$A$ with $\partial A\not\subset E$ we have to prove that $f$ is CR
near $\partial A$.  By symmetry it is enough to argue near $p=A(1)$.
We shall embed $A$ into family of analytic discs sweeping out a wedge
$\cal W$ glued to $M$. Using the continuity principle we shall extend
$f$ holomorphically to a subset of $\cal W$ of full measure and prove
that the extension $F$ is integrable on $\cal W$.

We shall first develope $A$ in a preliminary family $A_{\rho, s'}$ of
analytic discs, $0 \leq \rho < \rho_2$, $\rho_2 >\rho_1$,
$I_{\rho_2}=(0,\rho_2)$, $s'=(a_2,...,a_m,y_1^0,...,y_n^0)$ running in
a neighborhood of $0$ in $\C^{m-1}\times \R^m$, as follows: Set for
the $w$-component
\[
  W_{\rho,a}(\zeta)=(\rho(\zeta-\rho_1),a_2,...,a_{m-1})
\]
and take then $A_{\rho,s'}$ of the form
\[
  A_{\rho,s'}(\zeta)=(X_{\rho,s'}(\zeta)+
  iY_{\rho,s'}(\zeta),\rho(\zeta-\rho_1),a_2,...,a_{m-1}),
\]
where $Y_{\rho,s'}$ is the solution of Bishop's equation
\[
  Y_{\rho,s'}=T_1 h(Y_{\rho,s'},W_{\rho,s'})+y^0.
\]
Here $T_1$ denotes the renormalized Hilbert-transform which associates
to a (vector valued) real function on $\partial\Delta$ its harmonic
conjugate vanishing at the point $1\in\partial\Delta$. Since the discs
$W_{\rho,s'})$ are of small $C^{2,\alpha}$ norm, general properties of
Bishop's equation give the existence of solutions depending smoothly
on the data.

Differentiating Bishop's equation we see that there exists ${\cal V}$,
a neighborhood of $0$ in $\C^{m-1}$, ${\cal Y}$, a neighborhood of $0$
in $\R^n$, such that the mapping $$I_{\rho_2} \times {\cal A} \times
{\cal Y} \times b\Delta \ni (\rho,a,y^0,\zeta) \mapsto
A_{\rho,a,y^0}(\zeta) \in M$$ is an embedding. This shows that a
neighborhood in $M$ of $A_{0,0,0}(\overline{\Delta})=p_1$ is foliated
by ${\cal C}^2$-smooth real discs $D_{a,y^0}=D_{s'}=
\{A_{\rho,a,y^0}(\zeta)\in M\ \! {\bf :} \ \! 0\leq \rho < \rho_2,
\zeta \in b \Delta\}$.  Moreover, since $H_{d-2}(E)=0$, the set ${\cal
S}_E'=\{s'\in {\cal S}'\ \! {\bf :} \ \! D_{s'} \cap E \neq
\emptyset\}$ is a closed subset of ${\cal S}'={\cal A} \times {\cal
Y}$ of Lebesgue measure zero. By construction, each disc $A_{\rho,s'}$
with $s'\not\in {\cal S}_{E'}$ is therefore analytically isotopic to a
point in $M\backslash E$.

Furthermore, by means of normal deformations of the family near
$A(-1)$ as in [MP1], Proposition 2.6, we can develope $A$ in a regular
family $A_{\rho,s',v}$ which has the property that, for each $v\in
{\cal V}$, the set ${\cal S}_{E,v}'= \{s'\in {\cal S}'\ \! {\bf :} \
\! D_{s',v} \cap E\neq \emptyset\}$ is a closed subset of ${\cal S}'$
of Lebesgue measure zero.  Therefore, each disc $A_{\rho,s',v}$ with
$A_{\rho,s',v}(b\Delta)\cap E = \emptyset$ is analytically isotopic to
a point in $M\backslash E$, since ${\cal S}'\backslash {\cal
S}_{E,v}'$ is dense and open in ${\cal S}'$.

Then the isotopy property and a version of the continuity principle
(s.~[ME]) imply that ${\cal H}({\cal V}(M\backslash E))$ extends
holomorphically into
$${\cal W}=\{A_{\rho,a,y^0,v}(\zeta)\in \C^n\ \! {\bf :} \ \! \rho \in
I_{\rho_1}, a\in {\cal A}_1, y^0\in {\cal Y}_1, v\in {\cal V}_1, \zeta
\in \Delta_1\}$$ minus the set $$E_{{\cal W}}=\{A_{\rho,s',v}(\zeta)
\in \C^n\ \! {\bf :} \ \! A_{\rho,s',v}(b\Delta) \cap E \neq
\emptyset\}$$ for which $H_{2m+2n-1}(E_{{\cal W}}) = 0$.
% If even
% $H_{2m+2n-3}(E_{{\cal W}}) <\infty$, then $E_{{\cal W}}$ is removable
% and we are done.
Let $f\in {\cal H}(\omega)$ and let $F$ denote its extension to ${\cal
W}\backslash E_{{\cal W}}$.

Let us finally show $F\in L^1({\cal W})$. For each fixed $v$ the
boundaries of the discs $A_{\rho,s',v}$ foliate some open subset of
$M$. By applying a standard estimate to almost every of these discs
and integrating over $\rho,s'$ we get an $L^1$ estimate for the
restriction of $F$ to the $(d+1)$-dimensional submanifold of $W$ swept
out by these discs. Integrating over $v$ we finally get $F\in
L^1({\cal W})$.

{\it Step four: extension of $F$.}  The singular set $E_{{\cal W}}$
contains the closed subset
\[
  E'_{{\cal W}} = \{A_{\rho,s',v}(\zeta) \in \C^n\ \! {\bf :} \ \!
     A_{\rho,s',v}(b\Delta) \subset E \}.
\]
The main problem now consists in removing $E_{{\cal W}}\backslash
E'_{{\cal W}}$ from $W\backslash E'_{{\cal W}}$.  This is a
consequence of Theorem \ref{hol}. Indeed, the foliation of
$W\backslash E'_{{\cal W}}$ being given by the holomorphic discs
$A_{\rho,s',v}$ we observe that $E_{{\cal W}}\backslash E'_{{\cal W}}$
satisfies $H^{2N-1}(E_{{\cal W}}\backslash E'_{{\cal W}})=0$ because
of $H^{d-2}(E)=0$. Finally we have to verify that $F$ extends
holomorphically through some point of each leaf of $E_{{\cal
W}}\backslash E'_{{\cal W}}$. But that is clear by definition of
$E_{{\cal W}}\backslash E'_{{\cal W}}$ and our initial assumption
$f\in L_{loc}^1(M)\cap {\cal H}({\cal V}(M\backslash E))$.

Assuming Theorem \ref{hol} we are left with $E'_{{\cal W}}$.  Its
definition and $H^{d-2}(E)=0$ yield $H^{2N-2}(E'_{{\cal W}})=0$.
Hence almost every complex line is disjoint from $E'_{{\cal W}}$ and
we can extend $F$ through $E'_{{\cal W}}$ by the ordinary continuity
principle.

{\it Step five: $L^1$ boundary values.}  Finally, we wish to recover
$f$ near $p$ as the $L^1$-limit of $F$. For $0<\delta<<1,0<r<<1$,
define the approach manifolds
\[
  M_{v,r}=\{ A_{\rho,s',v}(re^{i\theta}):|1-\rho|<\delta, s'\in{\cal
  S},|\theta|<\delta \}.
\]
In the last but one step we proved $F\in L^1({\cal W})$ by estimating
along almost every disc. Similarly, we prove now a uniform $L^1$ bound
for the restrictions of $F$ to the approach manifolds.

Namely the Embedding-Theorem of Carleson (s.~[JO]) yields a uniform
bound
\[
  \int_{-\delta}^\delta |F(A_{\rho,s',v}(re^{i\theta}))|\; d\theta < C
  || f\circ A_{\rho,s',v} ||_{H^1},
\]
valid for almost all discs $A_{\rho,s',v}$ with $\partial
A_{\rho,s',v}\cap E=\emptyset$ (as usual $||\cdot||_{H^1}$ denotes the
Hardy-space norm).  Integrating over $s'$ gives the desired estimate
\[
  \int_{M_{v,r}} |F| \; d\theta ds' < C || f ||_{L^1(M')},
\]
where $M'$ is a sufficiently large neighborhood of $p$ in $M$.

As in the usual theory of Hardy spaces (s.~for the theory in several
variables the standard reference [ST]) $F$ attains near $p$ a weak
boundary $F^*$, which is an integrable CR-function.  Obviously, $F^*$
coincides (almost everywhere) with $f$, and the proof is ready.

\section{$L^1$-Removability for holomorphic functions}

We shall now complete the proof of Theorem \ref{main} by showing
Theorem \ref{hol}.  But before a remark concerning the relation of the
two theorems is in order.

\begin{rem}\rm
The reader may have noticed that we need in the proof of Theorem
\ref{main} only a weaker version of Theorem \ref{hol} where even
$H_{2N-1}(L)=0$ holds true. Under this additional hypothesis some
special cases get much simpler:

If $\cal F$ in Theorem \ref{hol} is even a {\sl holomophic} foliation, 
there is
an elementary proof by the continuity principle, which moreover does
not use the integrability of $F$.  For dimensional reasons the same is
true for $N=2$ (this argument generalizes readily to differential
foliations by analytic {\sl hypersurfaces}).  For the general
situation and $N\geq 3$ we do not know if a proof working only with
holomorphic hulls exists.
\end{rem}

{\it Proof of Theorem \ref{hol}.}  We take inspiration from an
argument used by Henkin and Tumanov (s.~[HT], Lemma 6) to treat a
related question for continuous CR-functions on CR-manifolds foliated
by complex curves.

Let $E_{nr}$ be the complement of the maximal open subset of $D$ to
which $F$ extends holomorphically. We assume $E_{nr}\not=\emptyset$
and have to deduce a contradiction.  Our arguments are local and work
near any point of $p_0\in E_{nr}$ which is not an inner point of the
set $L_{p_0}\cap E_{nr}$ with respect of the leaf-topology of the leaf
$L_{p_0}$ through $p_0$.

Around $p_0$ we may choose a neighborhood $Q$ with the following
properties: Near $\overline{Q}$ there are coordinates $w=u+iv$ and
$r=(r_1,\ldots,r_{N-2})$, where $w$ is holomorphic and the $r_j$ of
class $C^2$, such that $Q$ is given as $\{
0<u<1,0<v<1,0<r_j<1,j=1,\ldots,N-2 \}$ and the foliation $\cal F$
corresponds to the level sets of the mapping $r$. Contracting $U$
around $p_0$ we may further assume that $F$ is holomorphic near the
bottom $\{ v=0,0\leq u\leq 1,0\leq r_j\leq 1,j=1,\ldots,N-2 \}$.
Further we may assume the existence of holomorphic coordinates
$z=(z_1,\ldots,z_N)$ near $\overline{Q}$.  It shall be convenient to
work with slightly smaller product domain $Q'=\{0<v<1\}\times B'$
where we get $B'$ by smoothing the edges of the bottom.  In the
following we will tacitly suppose appropriate contractions of $Q$ and
$Q'$ around $p_0$ which does not destroy the precedingly achieved
properties.

We choose near $\overline{Q}$ a basis
$\overline{L_1},\ldots,\overline{L_N}$ of complex anti-linear
vectorfields whose coefficients are $C^2$ with respect to $(w,r)$ such
that $L_1$ is tangent to $\cal F$. As observed in [HT], the
$\overline{L_1},\ldots,\overline{L_N}$ may be corrected such that all
brackets of the form $[\overline{L_1},\overline{L_s}]$ are tangent to
$\cal F$.

Next we take a subdomain $G\subset Q'$ containing $p_0$ which is the
region squeezed between the bottom $B'$ and a smooth hypersurface $M$
which cuts $\{v=0\}$ transversely along $b B'$ and is transverse to
the leaves of $\cal F$.  By Fubini's Theorem, after a slight
deformation of $M$ the restriction $F|M$ may be supposed to be
integrable with respect to $(2N-1)$-dimensional volume.  We shall show
that (after some additional modifications) $F|M$ is an integrable
CR-function.  Afterwards, the usual Hartogs-Bochner-Theorem gives a
holomorphic extension of $F$ to $G$, in contradiction to the choice of
$p_0$. Let us rename $B:= B'$, $Q:= Q'$.

For technical reasons we have to fix in advance a special
approximation of $F$. Let $\chi$ be a smooth non-negative, compactly
supported, rotation-invariant function of the holomorphic coordinates
$z$ with $\int\chi dm(z)=1$. Take a smooth function $\eta\geq 0$ whose
support is contained in a small neighborhood of $\overline{Q}$ and
which equals $1$ near $\overline{Q}$.  Setting
$\chi_\eps(z)=(1/\eps)^{2N}\chi(z/\eps)$, we define for sufficiently
small $\eps>0$
\[
  F_\eps = (\eta F)*\chi_\eps,
\]
where $*$ denotes convolution with respect to Lebegues-measure in
$z$. It is standard that $F_\eps$ approximates $F$ in $L^1(Q)$. From
the mean-value property of holomorphic functions and the
rotation-invariance of $\chi$ we further deduce that, for near any
given point $z\in Q\backslash E_{nr}$, $F_\eps$ will coincide with $F$
for $\eps$ sufficiently small. In particular, this is true near the
bottom $B$.

We extract a subsequence $\eps_k\searrow 0$ and claim that after a
slight deformation of $M$ we can assume that the restricitions of
$F_k=F_{\eps_k}$ to $M$ tend in $L^1(M)$ to $M$. Indeed, this is a
consequence of the following variant of Fubini's theorem: If we have a
series of functions converging in $L^1$ on a product set, then their
restriction to almost all slices will converge to the restriction of
the limit.

Now we can return to the main part of the proof.  Fix a point $q\in
M\cap E_{nr}$.  We have to show
\begin{equation}\label{cr1}
  \int_M F\wedge\overline{\partial}\phi = 0,
\end{equation}
for any smooth $(N,N-2)$-form $\phi$ such that $\mbox{supp }\phi\cap
M$ is contained in some small neighborhood of $q$ in $M$. 
By the very definition of the tangential CR complex, the right
side of (\ref{cr1}) is not changed if we add to $\phi$ an 
$(N,N-2)$-form contained in the differential ideal generated
by a defining function $\rho$ of $M$ and its derivative 
$\overline{\partial}\rho$ (s.~[BOGG], 8.1). As $M$ is
transverse to $\cal F$, we may restrict our
attention to more special $\phi$:
Let $\overline{\omega_j}$ be a complex anti-linear dual base of $
\overline{L}_j$. Then the coefficient of $\overline{\omega_1}$
in $\overline{\partial}\rho$ is non-vanishing, and 
it is enough to consider $\phi=\sum'\phi_J
d\overline{\omega}_J$ such that all coefficients $\phi_J$ where $J$
contains $1$ are zero (here $\sum'$ means summation over increasing
indices, $\phi_J$ are $(N,0)$-forms).

Fixing such a $\phi$, we wish to approximate its coefficients (with
respect to the base $\omega_j, \overline{\omega_j}$) by functions
which are holomorphic along the leaves. In order to apply known
techniques for approximation in the complex plane we use the
holomorphic $w$-coordinate and argue fiberwise. As the intersection
of $\mbox{supp}\phi$ with a leaf $L_w$ is contained in a short segment
$I_w$, we may approximate the coefficient functions on $L_w$ by taking
convolution integrals over $I_w$ with holomorphic kernels (similar as
for instance in the proof of the Approximation-Theorem of Baouendi and
Treves).  As the integrals depend smoothly on $r$ we get that this
sequence of approximating functions $\phi_j$ tends in $C^2$ to $\phi$.
In particular, $\phi_j \to 0$ on $\partial M$.

The fact that $\mbox{supp}\phi_j$ can no longer be assumed to be of
compact support seems to cause complications.  Nevertheless, it shall
be possible to establish
\begin{equation}\label{cr2}
  \int_{M} F\wedge\overline{\partial}\phi_j = \int_{\partial M}
  F\wedge\phi_j,
\end{equation}
which implies (\ref{cr1}) by going to the limit, since $\phi_j \to 0$
on $\partial M$.

Now we verify as in [HT] that $\overline{\partial}
F_k\wedge\overline{\partial}\phi_j$ is of the form
$\overline{L_1}F_k\wedge\alpha$ where $\alpha$ is a $C^2$ form
indepent of $k$ (but of course depending on $j$). Indeed, for a
selection
$\overline{L_1},\overline{L_{s_1}},\ldots,\overline{L_{s_{N-2}}}$ we
apply Cartan's formula
\begin{eqnarray*}
  \lefteqn{(\overline{L_1},\overline{L_{s_1}},
     \ldots,\overline{L_{s_{N-2}}})\vdash\overline{\partial}\phi_j =}
     \\ & & \overline{L_1}(\overline{L_{s_1}},
     \ldots,\overline{L_{s_{N-2}}})\vdash\phi_j + \sum
     (-1)^{k+1}\overline{L_{s_k}} (\overline{L_1},\ldots,
     \widehat{\overline{L_{s_{k}}}},\ldots)\vdash\phi_j \\ & & +\sum
     (-1)^{k+h+1}(\ldots,[\overline{L_{s_k}},
     \overline{L_{s_h}}],\ldots)\vdash\phi_j +\sum
     (-1)^{k+1}([\overline{L_1},
     \overline{L_{s_k}}],\ldots)\vdash\phi_j,
\end{eqnarray*}
where $\vdash$ denotes interiour multiplication of vectors and forms.
By the choice of $\phi$ and of the $L_s$ such that $ [ \overline{L}_1,
\overline{L}_s ]= a_s \overline{L}_1$, all the right-hand terms
vanish. Consequently, $(\overline{L_1},\overline{L_{s_1}},
\ldots,\overline{L_{s_{N-2}}})\vdash\overline{\partial}\phi_j =0$,
that is to say, if we write $\overline{\partial}\phi_j=\sum'\beta_J
d\overline{\omega_J}$, all $\beta_J$ with $J$ containing $1$ vanish
identically.  Therefore in the wedge product with $\overline{\partial}
F_k = \sum_{j=1}^{N}\overline{L_i} F_k d\overline{\omega_i}$ only the
term $\overline{L_1} F_k d\overline{\omega_1}$ survives.

Consequently we get
\begin{equation}\label{cr3}
  \int_{M} F_k\wedge\overline{\partial}\phi_j = \int_{B}
  F_k\wedge\overline{\partial}\phi_j + \int_G \overline{L_1}
  F_k\wedge\alpha.
\end{equation}

So it remains to prove the following.

\smallskip
\noindent
{\bf Claim.} {\it For fixed $j$ and for $k\rightarrow\infty$, the
integral $\int_G \overline{L}_1 F_k\wedge\alpha$ tends to zero.}

\smallskip
Since by applying Stokes theorem afterwards, we will obtain (3):
$$\int_M F\wedge \overline{\partial} \phi_j = \int_B F \wedge
\overline{\partial} \phi_j \stackrel{\partial M = \partial B}{=}
\int_{\partial M} F\wedge \phi_j \longrightarrow_{j\infty} 0. $$

{\it Proof.} Let $R$ be the closed set of all
$r^0=(r^0_1,\ldots,r^0_{N-2})$ such that the holomorphic curve
$\{r=r^0\}$ has non-void intersection with $E_{nr}$. By a result of
Whitney, there is a so called regularized distance function
$\delta_R(r)$ which is comparable with the euclidean distance function
$\mbox{dist}(r,R)$ (s.~[ST]) and satisfies: (i) $1/C {\rm dist} (r, R)
\leq \delta_R(r) \leq C {\rm dist}(r, R)$, $\forall \ r\in B$, $C\geq
1$; (ii) $\delta_R(r) \in C^0(B) \cap C^2(B\backslash F)$; (iii)
$|\nabla \delta_R(r) | \leq C$ on $B\backslash F$.

Fix an arbitrary $\eps>0$. Further we take a small $\delta>0$ whose
precise choice shall be explained later. For the moment we only
require $\delta$ to be a regular value of $\delta_R$, i.~e.~to be
generic in the sense of the theorem of Sard. Hence
$Q\cap\{\delta_R=\delta\}$ a $C^2$-smooth hypersurface (fibered by
leaves of ${\cal F}$ and with possibly, a number of connected
components tending to $\infty$ as $\delta \to 0$.

Depending on $\delta$ we choose a large $k$ such that, for any $p\in
Q$ with $\delta_R(z)\geq\delta$, the distance of $p$ to the
singularity $E_{nr}$, measured with respect to the holomorphic
coordinate $z$, is greater than $\eps_k$ (the scaling parameter
appearing in the definition of $F_k$).

In the decomposition
\[
  \int_G \overline{L}_1 F_k\wedge\alpha =
  \int_{G\cap\{\delta_R>\delta\}} \overline{L}_1 F_k\wedge\alpha+
  \int_{G\cap\{\delta_R<\delta\}} \overline{L}_1 F_k\wedge\alpha
\]
the first term on the right is obviously zero since $F$ is holomorphic
on $G\cap\{\delta_R>\delta\}$. For the second term we integrate by
parts
\[
  \int_{G\cap\{\delta_R<\delta\}} \overline{L}_1 F_k\wedge\alpha =
  \int_{G\cap\{\delta_R<\delta\}} F_k\wedge\overline{L}_1^t\alpha +
  \int_{\partial(G\cap\{\delta_R<\delta\})} \sigma(\overline{L}_1)
  F_k\wedge\alpha,
\]
where $\overline{L}_1^t$ denotes the formally adjoint operator of
$\overline{L}_1$ and $\sigma(\overline{L}_1)$ the factor from the
symbol of $\overline{L}_1$ appearing in the boundary integral.

We claim that the first summand on the right gets small if $\delta$ is
small and $k$ large.  Indeed, as $R$ is of $(2N-2)$-dimensional volume
zero and $F$ integrable, there is a $\delta$ such that
\[
  \int_{Q\cap\{\delta_R<\delta\} } |F| < \eps.
\]
By the theorem of Sard, we may suppose that $\delta$ is a regular
value of $\delta_R$, and therefore $Q\cap\{\delta_R=\delta\}$ a smooth
hypersurface.  Next we choose $k$ so large that
\[
  \int_{Q\cap\{\delta_R<\delta\} } |F_k| < 2\eps.
\]

As $|\overline{L}_1^t\alpha| < C$ on $\overline{Q}$, we have
$|\int_{G\cap \{\delta_R < \delta\}} F_k \wedge \overline{L}_1^t
\alpha | \leq 2C \eps$. The second term decomposes in three boundary
integrals
\[
\int_{\partial (G \cap \{\delta_R < \delta\})} = \int_{G\cap
\{\delta_R=r\}} + \int_{M \cap \{\delta_R = r\}} + \int_{B\cap
\{\delta_R = r\}}
\]
For the first boundary term we remark that
\[
  \int_{G\cap\{\delta_R=\delta\}} \sigma(\overline{L}_1)
  F_k\wedge\alpha = 0,
\]
as $\overline{L}_1$ is tangential to $G\cap\{\delta_R=\delta\}$,
whence $\sigma(\overline{L}_1)$ vanishes on $G\cap \{\delta_R=r\}$.
The rest of the boundary integral is estimated in analogous manner as
the interiour integral. We have only to use that $F_k|M\rightarrow
F|M$ in $L^1$ and that $F_k$ coincides with $F$ near $B$. After
summing up the proof of the claim is finished.

\end{document}